\documentclass[reqno,12pt]{amsart}
\usepackage{geometry}
\usepackage{amsmath,amssymb,mathrsfs,color,stackengine}
\usepackage[colorlinks,
linkcolor=red,
anchorcolor=green,
citecolor=blue,
]{hyperref}
\usepackage[upint]{stix}
\usepackage{subfigure}
\usepackage{float}
\usepackage{times}
\usepackage{tikz}
\usetikzlibrary{intersections}
\usepackage{paralist}
\usepackage{comment}

\makeatletter
\def\setliststart#1{\setcounter{\@listctr}{#1}%
  \addtocounter{\@listctr}{-1}}
\makeatother

\makeatletter
\@addtoreset{figure}{section}
\makeatother

\usepackage{calc}
 \newtheorem{The}{Theorem}[section]
 \newtheorem{Cor}[The]{Corollary}
 \newtheorem{Lem}[The]{Lemma}
 \newtheorem{Pro}[The]{Proposition}
 \theoremstyle{definition}
 \newtheorem{defn}[The]{Definition}
 \newtheorem{Rem}[The]{Remark}
 \newtheorem{ex}[The]{Example}
 \numberwithin{equation}{section}

\newcommand{\R}{\mathbb{R}}

\newcommand{\N}{\mathbb{N}}

\title{topology of singular set of semiconcave function via Arnaud's theorem}
\author{Tianqi Shi, Wei Cheng \and Jiahui Hong}
\address{Department of Mathematics, Nanjing University, Nanjing 210093, China}
\email{tqshi.math@gmail.com}
\address{Department of Mathematics, Nanjing University, Nanjing 210093, China}
\email{chengwei@nju.edu.cn}
\address{Department of Mathematics, Nanjing University, Nanjing 210093, China}
\email{Jiahui.Hong.NJU@gmail.com}
\begin{document}
\maketitle

\begin{abstract}
We proved the (local) path-connectedness of certain subset of the singular set of semiconcave functions with linear modulus in general. In some sense this result is optimal. The proof is based on a theorem by Marie-Claude Arnaud (M.-C. Arnaud, \textit{Pseudographs and the Lax-Oleinik semi-group: a geometric and dynamical interpretation}. Nonlinearity, \textbf{24}(1): 71-78, 2011.). We also gave a new proof of the theorem in time-dependent case.
\end{abstract}

\section{Introduction}

The problem of singularities of semiconcave functions plays an important role in calculus of variation and optimal control, PDE, Hamiltonian dynamical systems and geometry. For convex Lagrangian on Euclidean space or finite dimensional manifold, the maximal regularity of the value function of the associated variational problem is semiconcavity. The study of the singularities of semiconcave functions in the literature can go back to the case of convex function which is a semiconvex function with linear modulus. The main theme of the previous study is the measure theoretic aspect of the set of singular points (\cite{Zajicek1979,Alberti_Ambrosio_Cannarsa1992,Albano_Cannarsa1999}) and the dynamic aspect, i.e., the propagation of singularities of semiconcave functions (\cite{Albano_Cannarsa2000,Cannarsa_Yu2009}). Recently, an intrinsic approach leads to some topological results of the cut locus as well as singular set of the viscosity solutions (\cite{Cannarsa_Cheng3,Cannarsa_Cheng_Fathi2017,Cannarsa_Cheng_Fathi2021}).

Let $M$ be a compact and connected smooth manifold without boundary. $TM$ and $T^*M$ are the tangent and cotangent bundle of $M$ respectively. A function $L=L(t,x,v):\R\times TM\to\R$ is called time-dependent \emph{Tonelli Lagrangian}, if it is a function of class $C^2$ and satisfies the following conditions:
\begin{enumerate}[(L1)]
	\item The function $v\mapsto L(t,x,v)$ is strictly convex for all $(t,x)\in\R\times M$.
	\item There exists a superlinear function $\theta:[0,+\infty)\to[0,+\infty)$ and $L^{\infty}_{\rm loc}$ function $c_0:\R\to[0,+\infty]$, such that
	\begin{align*}
		L(t,x,v)\geqslant \theta(|v|_x)-c_0,\qquad \forall (t,x,v)\in\R\times TM.
	\end{align*}
	\item There exist $C_1,C_2>0$, such that
	\begin{align*}
		|L_t(t,x,v)|\leqslant C_1+C_2L(t,x,v),\qquad \forall (t,x,v)\in\R\times TM.
	\end{align*}
\end{enumerate}
We call the Lagrangian $H$ associated with $L$ a \emph{Tonelli Hamiltonian}, defined by
\begin{align*}
	H(t,x,p)=\sup_{v\in T_xM}\{p(v)-L(t,x,v)\},\qquad (t,x,p)\in\R\times T^*M.
\end{align*}

Given any $x,y\in M$, and $t_1<t_2$, we denote by $\Gamma_{x,y}^{t_1,t_2}$ the set of absolutely continuous curves $\gamma\in \text{AC}([t_1,t_2],M)$ with $\gamma(t_1)=x$ and $\gamma(t_2)=y$. We call the function $h=h(t_1,t_2,x,y):\R\times\R\times M\times M\to\R$ the \emph{fundamental solution} of the associated Hamilton-Jacobi equation, define by
\begin{align*}
	h(t_1,t_2,x,y):=\inf_{\gamma\in\Gamma_{x,y}^{t_1,t_2}}\int_{t_1}^{t_2}L(s,\gamma(s),\dot\gamma(s))\,ds,\qquad t_1<t_2,\quad x,y\in M.
\end{align*}
For any $\phi:M\to\R$, $x\in M$, and $t_1<t_2$, we define
\begin{align*}
	T_{t_1}^{t_2}\phi(x):=\inf_{y\in M}\{\phi(y)+h(t_1,t_2,y,x)\},\\
	\breve T_{t_1}^{t_2}\phi(x):=\sup_{y\in M}\{\phi(y)-h(t_1,t_2,x,y)\}.
\end{align*}
The families $\{T_{t_1}^{t_2}\}$ and $\{\breve T_{t_1}^{t_2}\}$ of operators are called \emph{negative and positive Lax-Oleinik evolution} respectively. 

%\begin{Rem}
%If $\phi\in\text{\rm BUC}\,(M)$, the set of uniformly continuous functions on $M$, both $\{T_{t_1}^{t_2}\}$ and $\{\breve T_{t_1}^{t_2}\}$ satisfy Markov property. Moreover, if $L$ is independent of $t$, then both of them are semigroups on $\text{\rm BUC}\,(M)$ (see \cite{CCJWY2020}).
%\end{Rem} 

In this note, we will study the singular set of general semiconcave functions in the sense of topology, using some intrinsic observation based on a theorem by Marie-Claude Arnaud. We denote by $\Phi_H^{s,t}$ the associated flow of time-independent Hamiltonian $H$ from $s$ to $t$. Arnaud's theorem states that for $\phi\in\text{SCL}\,(M)$, the space of semiconcave functions (with linear modulus) on $M$ and $0<t_2-t_1\ll1$ such that $\breve{T}_{t_1}^{t_2}\phi\in C^{1,1}$, then
\begin{align*}
	\Phi_{H}^{t_1,t_2}(\text{\rm graph}\,(D\breve{T}_{t_1}^{t_2}\phi))=\text{\rm graph}\,(D^+\phi),
\end{align*}
where
\begin{align*}
	\text{graph}\,(D^+\phi)=&\,\{(x,p)\in T^*M: x\in M, p\in D^+\phi(x)\},\\
	\text{graph}\,(D\breve{T}_{t_1}^{t_2}\phi)=&\,\{(x,D\breve{T}_{t_1}^{t_2}\phi(x))\in T^*M: x\in M\}.
\end{align*}
In other words, there is a bi-Lipschitz diffeomorphism between the Lipschitz graph of $D\breve{T}_{t_1}^{t_2}\phi$ and the 1-graph $\text{graph}\,(D^+\phi)$. Thus, an immediate consequence is that $\text{graph}\,(D^+\phi)$ is a $d$-rectifiable set, which is also an $(\mathcal H^d,d)$-rectifiable set, $\Sigma^{\geqslant k}(\phi)=\{x\in M|\dim D^+\phi(x)\geqslant k\}$ is countably $(\mathcal H^{d-k},d-k)$-rectifiable, then so is $\Sigma^k(\phi)=\{x\in M|\dim D^+\phi(x)=k\}$. One can also have some easy topology consequences,
\begin{itemize}[$\bullet$]
	\item $\text{ind}\,(\Sigma^{\geqslant k}(\phi))\leqslant d-k$ for $k=0,\cdots,d$, where $\text{ind}\,(X)$ stands for the small inductive dimension of a $T_3$ topological space $X$.
	\item $\Sigma^{\geqslant k}(\phi)$ is an $F_\sigma$-set for $k=0,\cdots,d$. In particular, $\text{Sing}\,(\phi)=\Sigma^{\geqslant 1}(\phi)$ is an $F_\sigma$-set,   and in this case, $\text{Sing}\,(\phi)$ is of first category in the sense of Baire. $\Sigma^d(\phi)$ is at most countable.
\end{itemize}

In the case that a semiconcave function $u$ is a viscosity solution of Hamilton-Jacobi equation, $\text{Sing}\,(u)$ is even locally path-connected (\cite{Cannarsa_Cheng_Fathi2017,Cannarsa_Cheng_Fathi2021}). Main result of this paper is: Assume that $\phi\in\text{SCL}\,(M)$ and $d\geqslant 2$. Then for any $A\subseteq\Sigma^{\geqslant 2}(\phi)$, $M\setminus A$ is path-connected and locally path-connected. In particular, for $k=1,2,\cdots, d$, $\Sigma^{\leqslant k}(\phi)$ are path-connected and locally path-connected subsets of $M$. We also give an example to conclude that, in some sense, $\Sigma^{\leqslant 1}(\phi)=\Sigma^0(\phi)\cup\Sigma^1(\phi)$ is the smallest subset with path-connectedness and local path-connectedness, in general.

We also give a new proof of Arnaud's theorem for time-dependent Hamiltonian. This new proof also provides some new information on the intrinsic nature of the propagation of singularities and its connection to the non-commutativity of Lax-Oleinik commutators $T_{t_1}^{t_2}\circ\breve{T}_{t_1}^{t_2}$.

The paper is organized as follows. In section 2, we introduced some preliminary materials on semiconcave function, Lax-Oleinik operator and some basic notions on rectifiability and dimension theory from topology. In section 3, we proved the path-connectedness and other topological properties of $\text{Sing}\,(\phi)$ and subsets of $\text{Sing}\,(\phi)$ for $\phi\in\text{SCL}\,(M)$. We also gave a new proof of Arnaud's theorem for time-dependent Hamiltonian.

\medskip

\noindent\textbf{Acknowledgements.} Wei Cheng is partly supported by National Natural Science Foundation of China (Grant No. 12231010). %The authors also appreciate Kai Zhao for helpful discussion on the commutators.

\section{Preliminaries}

\subsection{Semiconcave functions}

We first recall some basic relevant notions on semiconcavity.
\begin{enumerate}[--]
	\item Let $\Omega$ be an open convex subset of $\R^d$. A function $\phi:\Omega\to\R$ is called a semiconcave function (of linear modulus) with constant $C\geqslant0$ if for any $x\in\Omega$ there exists $p\in\R^d$ such that
	\begin{equation}\label{eq:semiconcave}
	    \phi(y)\leqslant\phi(x)+\langle p,y-x\rangle+\frac C2|x-y|^2,\qquad \forall y\in\R^d.
	\end{equation}
	\item The set of covectors $p$ satisfying \eqref{eq:semiconcave} is called the proximal superdifferential of $\phi$ at $x$ and we denote it by $D^+\phi(x)$.
	\item Similarly, $\phi$ is semiconvex if $-\phi$ is semiconcave. The set of $D^-\phi(x)=-D^+(-\phi)(x)$ is called the proximal subdifferential of $\phi$ at $x$. 
	\item The set $D^+\phi(x)$ is a singleton if and only if $\phi$ is differentiable at $x$, and $D^+\phi(x)=\{D\phi(x)\}$. A point $x$ is called a singular point of a semiconcave function $\phi$ if $D^+\phi(x)$ is not a singleton. We denote by $\text{Sing}\,(\phi)$ the set of all singular points of $\phi$.
	\item We call $p\in D^*\phi(x)$, the set of reachable gradients, if there exists a sequence $x_k\to x$ as $k\to\infty$, $\phi$ is differentiable at each $x_k$ and $p=\lim_{k\to\infty}D\phi(x_k)$. We have $D^*\phi(x)\subset D^+\phi(x)$ and $D^+\phi(x)=\text{co}\,D^*\phi(x)$. 
\end{enumerate}
For more on the semiconcave functions, the readers can refer to \cite{Cannarsa_Sinestrari_book}.

\begin{Pro}\label{pro:inf}
A function $\phi:\Omega\to\R$ is a semiconcave function with constant $C\geqslant0$ if and only if there exists a family of $C^2$-function $\{\phi_i\}$ with $D^2\phi_i\leqslant CI$ such that
\begin{align*}
	\phi=\inf_i\phi_i.
\end{align*}
\end{Pro}

Proposition \ref{pro:inf} is an important and useful characterization of semiconcavity. Let $S$ be a compact topological space and $F:S\times\R^d\to\R$ be a continuous function such that
\begin{enumerate}[\rm (a)]
	\item $F(s,\cdot)$ is of class $C^2$ for all $s\in S$ and $\|F(s,\cdot)\|_{C^2}$ is uniformly bounded by some constant $C$,
	\item $D_xF(s,x)$ is continuous on $S\times\R^d$,
	\item $\phi(x)=\inf\{F(s,x): s\in S\}$.
\end{enumerate}
We call such a function $\phi$ a marginal function of a family of $C^2$-functions.

\begin{Pro}\label{pro:marginal}
Let $\phi(x)=\inf\{F(s,x): s\in S\}$ be a marginal function of the family of $C^2$-functions $\{F(s,\cdot)\}_{s\in S}$.
\begin{enumerate}[\rm (1)]
	\item $\phi$ is semiconcave with constant $C$.
	\item $D^*\phi(x)\subseteq\{D_xF(s,x): x\in\arg\min\{F(s,\cdot):s\in S\}\}$.
	\item For each $x$, let $M(x)=\arg\min\{F(s,x): s\in S\}$. Then
	\begin{align*}
	    D^+\phi(x)=
		\begin{cases}
			D_xF(s,x),& M(x)=\{s\}\ \text{is a singleton};\\
			\mbox{\rm co}\,\{D_xF(s,x): s\in M(x)\},&\text{otherwise.}
		\end{cases}
	\end{align*}
\end{enumerate} 	
\end{Pro}

\begin{Rem}
If $M$ is a connected and compact smooth manifold without boundary, a function $\phi:M\to\R$ is called semiconcave function, if there exists a family of $C^2$-functions $\{\phi_i\}$ such that $\phi=\inf_i\phi_i$, and the Hessians of $\phi_i$'s is uniformly bounded above. The readers can refer to the appendix of the paper \cite{Fathi_Figalli2010} for more  discussion of the semiconcavity of functions on (even noncompact) manifold. However, because of the local nature of our main discussion, we will work on Euclidean space instead. 
\end{Rem}

\subsection{Some aspects of rectifiability}

%\subsection{Some aspects of rectifiability}

For any metric space $(X,\rho)$ and $s\in[0,+\infty)$, let $\mathcal H^s$ be the $s$-dimensional Hausdorff measure in $X$. For any $E\subset X$,  The \emph{Hausdorff dimension} of a subset $E$ is
\begin{align*}
	\dim_H(E):=\inf\{s\geqslant 0|\mathcal H^s(E)=0\}=\sup\{s\geqslant 0|\mathcal H^s(E)=+\infty\}.
\end{align*}

%We first recall the definition of the $s$-dimensional Hausdorff measure in $\R^d$. For $E\subseteq\R^d$, $s\in[0,+\infty)$ and $\delta\in(0,+\infty]$, we set
%\begin{align*}
%	\mathcal H_\delta^s(E):=\left\{\sum_{i=1}^\infty\alpha(s)\left(\frac{\text{diam}C_i}{2}\right)^s\bigg|E\subseteq \bigcup_{i=1}^\infty C_i,\text{diam}C_i\leqslant\delta\right\},
%\end{align*}
%where $\text{diam}C_i:=\sup\{|x-y||x,y\in C_i\}$ and 
%\begin{align*}
%	\alpha(s)=\frac{\pi^{\frac{s}{2}}}{\Gamma\left(\frac{s}{2}+1\right)}
%\end{align*}
%and $\Gamma(s)=\int_0^{+\infty}e^{-x}x^{s-1}dx$ is the Gamma function. The following and further related definition and properties can be found in \cite{Evans_Gariepy_book}.
%
%\begin{defn}
%	($s$-dimensional Hausdorff measure) For any $E\subseteq \R^d$ and $s\in[0,+\infty)$, the $s$-\emph{dimensional Hausdorff measure} of $E$ is defined as
%	\begin{align*}
%	\mathcal H^s(E):=\lim_{\delta\to 0}\mathcal H_{\delta}^s(E)=\sup_{\delta>0}	\mathcal H_{\delta}^s(E).
%	\end{align*}
%\end{defn}
%
%\begin{Rem}
%	For a general metric space $(X,d(\cdot,\cdot))$ and $U\subseteq X$ is an nonempty set. Define the \emph{diameter} of $U$ as $\text{diam} U:=\sup\{d(x,y)|x,y\in U\}$. Then we set
%	\begin{align*}
%	\mathcal H_\delta^s(E):=\left\{\sum_{i=1}^\infty(\text{diam}C_i)^s\bigg|E\subseteq \bigcup_{i=1}^\infty C_i,\text{diam}C_i\leqslant\delta\right\}.
%\end{align*}Similarly, the $s$-dimensional Hausdorff measure in $(X,d(\cdot,\cdot))$ can be defined as above.
%\end{Rem}

\begin{Pro}\label{haus meas}
The Hausdorff measure satisfies the following properties.
\begin{enumerate}[\rm (1)]
	\item The measure $\mathcal H^s$ is Borel regular measure for all $s\in[0,+\infty)$;
	\item $\mathcal H^d$ coincides with the Lebesgue measure of $\R^d$ and $\mathcal H^0$ is the counting measure;
	\item If $f\in \emph{Lip}(\R^m;\R^d)$ and $E\subseteq\R^m$, then we have $\mathcal H^s(f(E))\leqslant(\emph{Lip}\,(f))^s\mathcal H^s(E)$;
	\item Let $E\subseteq\R^d$ and $0\leqslant s<t<+\infty$. If $\mathcal H^s(E)<+\infty$, then $\mathcal H^t(E)=0$. If $\mathcal H^t(E)>0$, then $\mathcal H^s(E)=+\infty$.
\end{enumerate}
\end{Pro}

%Thanks to the properties above, the definition below makes sense.

\begin{Pro}\label{prop of hausdorff dimension}\
\begin{enumerate}[\rm (1)]
	\item $\forall E\subseteq \R^d$, $\dim_HE\leqslant d$. And for $s=\dim_HE$, $\mathcal H^s(E)<+\infty$ or $\mathcal H^s(E)>0$;
	\item Assume that $F_i\subseteq\R^d$ for $i=1,2,\cdots$, then 
	\begin{align*} 
	\dim_H(\cup_{i=1}^\infty F_i)=\sup_{1\leqslant i\leqslant\infty}\{\dim_H(F_i)\}.
	\end{align*}
\end{enumerate}
\end{Pro}

The following definition of rectifibaility with respect to Hausdorff measure comes from \cite{Federer_book1969}.
\begin{defn}(Rectifiability)\label{rec}
	Assume that $E\subseteq\R^d$ and $m$ is a positive integer.
	\begin{enumerate}[\rm (1)]
\item $E$ is $m$-\emph{rectifiable} if and only if there exists a Lipschitz map $f:\R^m\to\R^d$ and a bounded set $B\subseteq \R^m$, such that $f(B)=E$;
\item $E$ is \emph{countably} $m$-\emph{rectifiable} if and only if $E$ equals the union of some countable family whose members are $m$-rectifiable;
\item $E$ is \emph{countably} $(\mathcal H^m,m)$-\emph{rectifiable} if and only if there exists a countably $m$-rectifiable set containing $\mathcal H^m$ almost all of $E$;
\item $E$ is $(\mathcal H^m,m)$-\emph{rectifiable} if and only if $E$ is countably $(\mathcal H^m,m)$-rectifiable and $\mathcal H^m(E)<+\infty$.
\end{enumerate}
\end{defn}

\begin{Rem}\label{rec} 
For rectifiablity, we have some comments as follows.
	\begin{enumerate}[(1)]
		\item If $E\subseteq\R^d$ is $m$-rectifiable, then $E$ is $(\mathcal H^m,m)$-rectifiable;
		
		This is straightforward from the facts that $E$ is the image of a bounded subset of $\R^m$ under a Lipschitz map $f$ and 
		\begin{align*}
		\mathcal H^m(E)=\mathcal H^m(f(B))\leqslant(\text{Lip}\,(f))^m\mathcal H^m(B)<+\infty.	
		\end{align*}
		\item If $E\subseteq\R^d$ is $(\mathcal H^m,m)$-rectifiable, then $\dim_H(E)\leqslant m$;
		\item $E\subseteq\R^d$ is countably $(\mathcal H^m,m)$-rectifiable if and only if $E$ can be represented as a countable union of $(\mathcal H^m,m)$-rectifiable sets;
		
		Since $E$ is countably $(\mathcal H^m,m)$-rectifiable, it is equivalent to there exists a sequence of $m$-rectifiable sets $V_i\subseteq\R^d$,
		\begin{align*}
		\mathcal H^m\left(E-\bigcup_{i=1}^\infty V_i\right)=0.	
		\end{align*}
 Now let $V_0=E-\cup_{i=1}^\infty V_i$, then $V_0$ is $(\mathcal H^m,m)$-rectifiable, and $E=\cup_{i=0}^\infty V_i$. Conversely, If $\{V_i\}_{i=1}^\infty$ is the sequence of the $(\mathcal H^m,m)$-rectifiable set and $E=\cup_{i=1}^\infty V_i$, then for each $i$, there exists $\{V_{ij}\}_{j=1}^\infty$ are the $m$-rectifiable sets, which satisfies
 \begin{align*}
 	\mathcal H^m\left(V_i-\bigcup_{j=1}^\infty V_{ij}\right)=0.	
 \end{align*}
Thus we can choose countable $m$-rectifiable sets $\{V_{ij}\}_{i,j=1}^{\infty}$, such that
 \begin{align*}
		\mathcal H^m\left(E-\bigcup_{i,j=1}^\infty V_{ij}\right)\leqslant\sum_{i=1}^\infty\mathcal H^m\left(V_i-\bigcup_{j=1}^\infty V_{ij}\right)=0.	
		\end{align*}
		\item If $E\subseteq\R^d$ is a (countably) $(\mathcal H^m,m)$-rectifiable set, $F\subseteq E$, then so is $F$.
	\end{enumerate}
\end{Rem}

\begin{Pro}\label{federer}
	If $W$ is an $(\mathcal H^m,m)$-rectifiable set and $\mathcal H^m$ measurable subset of $\R^k $, $f:W\to\R^\nu$ is a Lipschitz map, $\mu$ is an integer, $0\leqslant \mu\leqslant m$, $\lambda>0$ and 
	\begin{align*}
	Z_\lambda:=\{y\in\R^{\nu}|\mathcal H^{m-\mu}(f^{-1}(y))\geqslant\lambda\},
	\end{align*}
then $Z_\lambda$ is $(\mathcal H^\mu,\mu)$-rectifiable. Thus
\begin{align*}
	Z_0:=\{y\in\R^\nu|\mathcal H^{m-\mu}(f^{-1}(y))>0\}
	\end{align*}
	is countably $(\mathcal H^\mu,\mu)$-rectifiable. 
\end{Pro}

\begin{proof}
The proof of $(\mathcal H^\mu,\mu)$-rectifiablity of $Z_\lambda$ is due to Herbert Federer (See \cite[Theorem 3.2.31]{Federer_book1969}). To prove the rest of the argument, we only need to notice that $Z_0=\bigcup_{n=1}^\infty Z_{\frac{1}{n}}$ and (3) of Remark \ref{rec} implies the conclusion.
\end{proof}

\subsection{General topology and dimension theory}

Assume that $X$ is a topological space.
\begin{defn}
	(small inductive dimension, see \cite[Definition 1.1.1]{Engelking_book1978}) Assume that $X$ is a $T_3$ topological space, we define the \emph{small inductive dimension} denoted by $\text{ind}\,(X)$,  an integer "inductively":
	\begin{enumerate}[\rm (1)]
		\item $\text{ind}\,(X)=-1$ if and only if $X=\varnothing$;
		\item $\text{ind}\,(X)\leqslant n$, where $n\in\N^*$, if for every point $x\in X$ and each neighbourhood $V\subseteq X$ of the point $x$ there exists an open set $U\subseteq X$ such that $x\in U\subseteq V$ and $\text{ind}\,(\partial U)\leqslant n-1$;
		\item $\text{ind}\,(X)= n$ if $\text{ind}\,(X)\leqslant n$ and $\text{ind}\,(X)>n-1$, i.e. the inequality $\text{ind}\,(X)\leqslant n-1$ does not hold;
		\item $\text{ind}\,(X)=+\infty$ if $\text{ind}\,(X)>n$ for all $n=-1,0,1,\cdots$.
	\end{enumerate}
\end{defn}

%\begin{defn}
%	(large inductive dimension, see \cite{engelking_dimension_1978}, Definition 1.6.1) Assume that $X$ is a $T_4$ topological space, we define the \emph{large inductive dimension} denoted by $\text{Ind} X$,  an integer "inductively":
%	\begin{itemize}
%		\item[(1)] $\text{Ind} X=-1$ if and only if $X=\varnothing$;
%		\item[(2)] $\text{Ind} X\leqslant n$, where $n\in\N^*$, if for every closed set $A\subseteq X$ and each open set $V\subseteq X$ which contains $A$ there exists an open set $U\subseteq X$ such that $A\subseteq U\subseteq V$ and $\text{Ind}\,\partial U\leqslant n-1$;
%		\item[(3)] $\text{Ind} X= n$ if $\text{Ind} X\leqslant n$ and $\text{Ind} X>n-1$, i.e. the inequality $\text{Ind} X\leqslant n-1$ does not hold;
%		\item[(4)] $\text{Ind}X=+\infty$ if $\text{Ind} X>n$ for all $n=-1,0,1,\cdots$.
%\end{itemize}
%\end{defn}

\begin{The}\label{dimension}\
	\begin{enumerate}[\rm (1)]
		
		\item For any metric space $X$, $\emph{ind}\,(X)\leqslant \dim_H(X)$;
		\item For any $\varnothing\neq E\subseteq\R^d$, we have $0\leqslant \emph{ind}\,(E)\leqslant\dim_H(E)\leqslant d$;
		\item For every subspace $Y$ of a separable metric space $X$, we have $\emph{ind}\,(Y)\leqslant \emph{ind}\,(X)$.
	\end{enumerate}
\end{The}

\begin{proof}
	The conclusion (1) follows from  Theorem 6.3.10 of \cite{Edgar_book2008}. (2) can be obtained directly from (1). (3) is a direct corollary of Theorem 1.1.2 in \cite{Engelking_book1978}.
\end{proof}

\begin{Pro}\label{convex dim}
	For a convex set $E\subseteq\R^d$, $\dim\,(E)=\emph{ind}\,(E)=\dim_H(E)$, where $\dim\,(E)$ is the dimension of convex set $E$ (See \cite[Definition 2.1.1]{Hiriart-Urruty_Lemarechal_book2001}).
\end{Pro}

\begin{proof}
	For each point $x$ in the relative interior of $E$, which we denote by $\text{ri} E$ (See \cite[Section 2.1]{Hiriart-Urruty_Lemarechal_book2001} for more information), there exists $\delta>0$ such that the intersection of the affine hull $\text{aff}E$ and $d$-dimensional closed ball $B^d(x,\delta)$ contains in $E$. And $\text{ind}\,(\partial(\text{aff}E\cap B^d(x,\delta)))=\text{ind}\,(S^{\dim(\text{aff}E)-1})=\dim\,(\text{aff}E)-1$. The arbitrariness of $x\in\text{ri}E$ and Theorem \ref{dimension} imply that 
	\begin{align*}
		\text{ind}\,(E)\geqslant\text{ind}\,(\text{ri}E)=\dim\,(\text{aff}\,E)=\dim\,(E).
	\end{align*}
	On the other hand, 
	\begin{align*}
		\dim\,(E)=\dim\,(\text{aff}\,E)=\text{ind}\,(\text{aff}\,E)\geqslant \text{ind}\,(E).
	\end{align*}Consequently, we have $\dim\,(E)=\text{ind}\,(E)$. Moreover, with Theorem \ref{dimension},
	\begin{align*}
	\dim_H(E)\geqslant\text{ind}\,(E)=\dim\,(E)=\dim\,(\text{aff}\,E)=\dim_H(\text{aff}\,E)\geqslant \dim_H(E),	
	\end{align*}
	which completes the proof.
\end{proof}

\section{Topology of the singular set of semiconcave functions}

Let $M$ be a connected and closed manifold of dimension $d$. Let $\text{SCL}\,(M)$ be the set of semiconcave functions on $M$. Let $H=H(x,p)$ be any time-independent Tonelli Hamiltonian on $M$. This allows us to apply following Proposition \ref{pro:Arnaud} directly. If $t_1<t_2$, then $\breve{T}^{t_2}_{t_1}=\breve{T}^{t_2+s}_{t_1+s}$ for any $s\in\R$. Thus, we set $t=t_2-t_1$ and $\breve{T}^{t_2}_{t_1}=\breve{T}_0^t$. For any semiconcave function $\phi\in\text{SCL}\,(M)$, let 
\begin{align*}
	t_\phi:=&\,\sup\{t>0:  \breve{T}_0^t\phi\in C^{1,1}(M)\}.
\end{align*}
For any $\phi\in\text{SCL}(M)$, and $t\in(0,t_\phi]$, we denote by
\begin{align*}
	\text{graph}\,(D^+\phi)=&\,\{(x,p)\in T^*M: x\in M, p\in D^+\phi(x)\},\\
	\text{graph}\,(D\breve{T}_0^{t}\phi)=&\,\{(x,D\breve{T}_0^{t}\phi(x))\in T^*M: x\in M\}.
\end{align*}
We call the set $\text{\rm graph}\,(D^+\phi)$ the \emph{pseudo-graph} of $ D\phi$ or \emph{1-graph} of $\phi$. Now, let us recall a result by Marie-Claude Arnaud (\cite{Arnaud2011}) at first.

\begin{Pro}\label{pro:Arnaud}
Let $\phi\in\text{\rm SCL}\,(M)$ and $t_0>0$ such that $\breve{T}^{t_0}_{0}\phi\in C^{1,1}(M)$. Then, $\breve{T}_0^{t}\phi\in C^{1,1}(M)$ for all $t\in(0,t_0]$, and we have
	\begin{equation}\label{eq:graph_evo}
		\text{\rm graph}\,(D\breve{T}^t_0\phi)=\Phi_H^{-t}(\text{\rm graph}\,(D^+\phi)),\qquad\forall t\in(0,t_0],
	\end{equation}
	where $\{\Phi^t_H\}_{t\in\R}$ is the Hamiltonian flow with respect to any time-independent Tonelli Hamiltonian $H:T^*M\to\R$.
\end{Pro}

The original proof theorem of Arnaud is based some approximation of the graph for time-independent Lagrangian. We give a new proof of a refinement (Theorem \ref{thm:C11}) of Proposition \ref{pro:Arnaud} in time dependent case with more geometric intuition. The relation \eqref{eq:graph_evo} defines a bi-Lipschitz diffeomorphism between $\text{\rm graph}\,(D^+\phi)$ and $\text{\rm graph}\,(D\breve T_0^t\phi)$ by the Hamiltonian flow. This allows us to understand the singularities of $\phi$ by means of the latter Lipschitz graph. Expect for some finer results obtain in \cite{Albano_Cannarsa1999}, this idea also leads to more topological results of the singular set of semiconcave functions.  % and $c$-concave function from the theory of optimal transport.

%Before our discussion on some topological and measure theoretic aspects of the singular set of semiconcave functions, we begin with some problem from the issue of propagation of singularities.
%
%\begin{Lem}\label{argmin}
%	For any $\phi\in\emph{SCL}(M)$, $t\in\R$ and $x\in M$, 
%\begin{align*}
%	\lim_{\varepsilon\to 0^+}D\breve T_{t}^{t+\varepsilon}\phi(x)=\mathop{\arg\min}_{p\in D^+\phi(x)}H(t,x,p).
%\end{align*}
%\end{Lem}

\subsection{Singular set of semiconcave functions}

\begin{Lem}
	For any $\phi\in\emph{SCL}\,(M)$, $\emph{graph}\,(D^+\phi)$ is a $d$-rectifiable set, which is also an $(\mathcal H^d,d)$-rectifiable set. 
\end{Lem}
\begin{proof}
	Since for sufficiently $t\in(0,t_\phi]$, $\breve T_0^t\phi\in C^{1,1}(M)$, then by the Lipschitz map:
	\begin{align*}
		(\text{id}\times D\breve T_0^t\phi):M&\to\text{graph}\,(D\breve T_0^t\phi)\subseteq T^*M\\
		x&\mapsto(x,D\breve T_0^t\phi(x)),
	\end{align*} we can obtain that $\text{graph}\,(D\breve T_0^t\phi)$ is a $d$-rectifiable set in $T^*M$, because of the compactness of $M$. According to Theorem \ref{thm:C11}, for all $t_2-t_1\in(0,t_\phi]$, 
\begin{align*}
	\text{graph}\,(D^+\phi)=\Phi_{H}^{t}(\text{graph}\,(D\breve T_0^t\phi)),
\end{align*}and the Lipschitz property of the Hamiltonian flow $\Phi_{H}^t$ on $T^*M$, we know from the Definition \ref{rec} that $\text{graph}\,( D^+\phi)$ is a $d$-rectifiable set and thus an $(\mathcal H^d,d)$-rectifiable set.
\end{proof}

For any $\phi\in\text{SCL}\,(M)$, and integer $k$ with $0\leqslant k\leqslant d$, we define 
\begin{align*}
\Sigma^{\geqslant k}(\phi)&:=\{x\in M|\dim D^+\phi(x)\geqslant k\},\\
\Sigma^{\leqslant k}(\phi)&:=\{x\in M|\dim D^+\phi(x)\leqslant k\},\\
\Sigma^{k}(\phi)&:=\{x\in M|\dim D^+\phi(x)= k\}.
\end{align*}
 Obviously, $\Sigma^0(\phi)$ is the set of differentiable points of $\phi$, $\Sigma^{\geqslant 0}(\phi)=\Sigma^{\leqslant d}(\phi)=M$ and $\Sigma^{\geqslant 1}(\phi)=\text{Sing}\,(\phi)$.

\begin{Pro}\label{d-k}
	$\Sigma^{\geqslant k}(\phi)$ is countably $(\mathcal H^{d-k},d-k)$-rectifiable for $k=0,\cdots,d$, then so is $\Sigma^k(\phi)$.
\end{Pro}

\begin{proof}
	 Consider the canonical projection map
	\begin{align*}
		\pi_x:\text{graph}(D^+\phi)&\to M\subseteq\R^{2d}\\
		(x,p)&\mapsto x,
	\end{align*}
	which is a Lipschitz map. In this case, we have $\pi_x^{-1}(x)=(x, D^+\phi(x))$, for all $x\in M$. Thus, by Theorem \ref{federer}, for any $k=0,1,\cdots,d$, $S_k:=\{x\in M|\mathcal H^{d-(d-k)}(\pi_x^{-1}(x))>0\}$ is a countably $(\mathcal H^{d-k},d-k)$-rectifiable set. Moreover from Proposition \ref{convex dim}, we know that
		\begin{align*}
		\Sigma^{\geqslant k}(\phi)&=\{x\in M|\dim\,(D^+\phi(x))\geqslant k\}\\
				&=\{x\in M|\dim_H(D^+\phi(x))\geqslant k\}\\
				&=\{x\in M|\mathcal H^{k}((x,D^+\phi(x)))>0\}\\
				&=\{x\in M|\mathcal H^{d-(d-k)}(\pi_x^{-1}(x))>0\}\\
				&=S_k.
			\end{align*}Consequently, $\Sigma^{\geqslant k}(\phi)$ is a countably $(\mathcal H^{d-k},d-k)$-rectifiable set, thus so is $\Sigma^{k}(\phi)$ by (4) of Remark \ref{rec}.
\end{proof}

\begin{Rem}\label{tangent cone}
Actually, one can show that $\Sigma^k(\phi)$ is a countably $(d-k)$-rectifiable set in Euclidean case by analyzing the tangent cone of $\Sigma^k(\phi)$. For further reading, the reader can refer to \cite[Section 4.1]{Cannarsa_Sinestrari_book}.
\end{Rem}

\begin{Cor}
	Assume that $\phi\in\emph{SCL}\,(M)$, then $\emph{ind}\,(\Sigma^{\geqslant k}(\phi))\leqslant d-k$ for $k=0,\cdots,d$.
\end{Cor}

\begin{proof}
From Theorem \ref{d-k}, we know that $\Sigma^{\geqslant k}(\phi)$ is a countably $(\mathcal H^{d-k},d-k)$-rectifiable set, then there exists a sequence of $(\mathcal H^{d-k},d-k)$-rectifiable sets $\{F_i\}_{i=1}^\infty$, i.e., $\dim_H(F_i)\leqslant d-k$, such that $\Sigma^{\geqslant k}(\phi)=\bigcup_{i=1}^\infty F_i$. Proposition \ref{prop of hausdorff dimension} and Theorem \ref{dimension} imply that
\begin{align*}
\text{ind}\,(\Sigma^{\geqslant k}(\phi))\leqslant\dim_H(\Sigma^{\geqslant k}(\phi))=\sup_{1\leqslant i\leqslant\infty}\{\dim_H (F_i)\}\leqslant d-k.	
\end{align*}
This completes the proof.
\end{proof}

\begin{Rem}\label{connected}
It is well known that when $d\geqslant 2$, for any $d$-dimensional connect and compact manifold $M$ and subset $A$ with $\text{ind}\,(A)\leqslant {d-2}$, $M\setminus A$ is connected (See \cite[Example VI.11]{Hurewicz_Wallman1941}). So from the corollary above, we are clear that if $A\subseteq\Sigma^{\geqslant 2}(\phi)$, $M\setminus A$ is connected, especially $\Sigma^{\leqslant k}(\phi)$ are also connected subsets for $k=1,2,\cdots,d$.
\end{Rem}

In addition to dimension, we also have further topological conclusions about the singularities of semiconcave functions.

\begin{Pro}\label{fsigma_sing}
	For any $\phi\in\emph{SCL}\,(M)$, $k=0,1,\cdots, d$, $\Sigma^{\geqslant k}(\phi)$ is an $F_\sigma$-set. In particular, $\emph{Sing}(\phi)=\Sigma^{\geqslant 1}(\phi)$ is an $F_\sigma$-set,   and in this case, $\emph{Sing}\,(\phi)$ is of first category in the sense of Baire. 
\end{Pro}

\begin{proof}
This is a direct consequence of Proposition 4.1.11 of \cite{Cannarsa_Sinestrari_book} and Rademencher's theorem.
\end{proof}

\begin{Pro}\label{countable}
Let $\phi\in\emph{SCL}\,(M)$, then $\Sigma^d(\phi)$ is at most countable, thus totally disconnected when it is not empty, i.e. the connected components of $\Sigma^d(\phi)$ are singletons.
\end{Pro}

\begin{proof}
According to the definition of $\Sigma^d(\phi)$, we can know that for each $x\in\Sigma^d(\phi)$ and $t\in(0,t_\phi]$, by Theorem \ref{thm:C11}, $\Phi_{H}^{-t}$ is a $C^1$-diffeomorphism, then
\begin{align*}
	K_x:=\Phi_{H}^{-t}(x,D^+\phi(x))\subseteq\text{graph}(D\breve T_{0}^{t}\phi)
\end{align*}is a $d$-dimensional closed region. And for two distinct $x,y$ in $\Sigma^d(\phi)$, $K_x\cap K_y=\varnothing$, which means $\Sigma^d(\phi)$ is at most countable. Thus, $\Sigma^d(\phi)$ is totally disconnected.
\end{proof}

\begin{The}\label{path}
Assume that $\phi\in\emph{SCL}\,(M)$ and $d\geqslant 2$. For any distinct $a,b\in M$, there exists a piecewise smooth (broken line) curve $\gamma: [0,1]\to M$ with $\gamma(0)=a$, $\gamma(1)=b$ and $\gamma(s)\notin\Sigma^{\geqslant 2}(\phi)$ for $s\in(0,1)$.
\end{The}

\begin{proof}
Since $M$ is a connected manifold, we can discuss in a local coordinate neighborhood to complete the proof of the theorem by splicing curves. For convenience, we can assume that $M=B^d_R$, the $d$-dimensional closed ball with radius $R>0$, for the case of dimension $d$.

Now for two distinct points $a,b\in B_R^d$,  we can find a $d$-dimensional closed disk $D_r^{d-1}\subseteq B^d_R$ with $r\in(0,R)$, such that $a$ and $b$ are on opposite sides of $\text{aff}\,(D_r^{d-1})$. For any fixed $z\in D_r^{d-1}$, we define a broken line $l_z:[0,1]\to B_R^d$ as
\begin{align*}
	l_z(s)=
	\begin{cases}
		(1-2s)a+2sz,\,\,\,\,\,\,\,\,\,\,\,\,\;s\in[0,\frac{1}{2}],\\
		(2-2s)z+(2s-1)b, s\in[\frac{1}{2},1].
	\end{cases}
\end{align*}
Set the map
\begin{align*}
	S:C_{a,b}^d&\to D_{r}^{d-1}\\
	l_z\ni x&\mapsto z,
\end{align*}
where $C_{a,b}^d=\{l_z(s)|s\in(0,1),z\in D_r^{d-1}\}$. $S$ is well defined since $l_{z_1}\cap l_{z_2}=\{a,b\}$ for any $z_1\neq z_2$. 
Let $1\gg r_n\searrow 0$ and $C_n:=C_{a,b}^d\setminus (B^d(a,r_n)\cup B^d(b,r_n))$, then for any $n\in\N^*$, $C_n$ is compact and $C_{a,b}^d=\cup_{n=1}^\infty C_n$. 

Recall that $\Sigma^{\geqslant 2}(\phi)$ is a countably $(\mathcal H^{d-2},d-2)$-rectifiable set by Proposition \ref{d-k}.  Then  combining with (3) in Remark \ref{rec}, we can find a sequence of $(\mathcal H^{d-2},d-2)$-rectifiable sets $\{W_m\}_{m=1}^\infty$, such that $\Sigma^{\geqslant 2}(\phi)=\cup_{m=1}^\infty W_m$. Thus
\begin{align*}
	S\left(\Sigma^{\geqslant 2}(\phi)\cap C_{a,b}^d\right)=S\left(\bigcup_{n=1}^\infty\bigcup_{m=1}^\infty(W_m\cap C_n)\right)=\bigcup_{n=1}^\infty\bigcup_{m=1}^\infty S(W_m\cap C_n).
\end{align*}Note that $W_m\cap C_n$ is an $(\mathcal H^{d-2},d-2)$-rectifiable set for each $m,n\in\N^*$, which implies that $\mathcal H^{d-2}(W_m\cap C_n)<+\infty$ and  $\mathcal H^{d-1}(W_m\cap C_n)=0$ by (4) of Proposition \ref{haus meas}. It is easy to check that for any fixed $n$, $S|_{C_n}$ is a Lipschitz map. Thus, Proposition \ref{fsigma_sing} implies that $S(\Sigma^{\geqslant 2}(\phi)\cap C_{a,b}^d)$ is an $F_\sigma$-set. 
Moreover, from (3) of Proposition \ref{haus meas}, we can infer that
\begin{align*}
	\mathcal H^{d-1}(S(W_m\cap C_n))\leqslant\left(\text{Lip}\left(S|_{c_n}\right)\right)^{d-1}\mathcal H^{d-1}(W_m\cap C_n)=0.
\end{align*}According to (2) in Proposition \ref{haus meas}, $\mathcal H^{d-1}$ is actually the $(d-1)$-dimensional Lebesgue measure $\mathscr L^{d-1}$ on $D_r^{d-1}$. Consequently, we have
\begin{align*}
	\mathscr L^{d-1}\left(S\left(\Sigma^{\geqslant 2}(\phi)\cap C_{a,b}^d\right)\right)&=\mathcal H^{d-1}\left(S\left(\Sigma^{\geqslant 2}(\phi)\cap C_{a,b}^d\right)\right)\\
	&=\mathcal H^{d-1}\left(\bigcup_{n=1}^\infty\bigcup_{m=1}^\infty S(W_m\cap C_n)\right)\\
	&\leqslant\sum_{n=1}^\infty\sum_{m=1}^\infty \mathcal H^{d-1}\left(S(W_m\cap C_n)\right)\leqslant0\\
	&<\mathcal H^{d-1}(D_r^{d-1})=\mathscr L^{d-1}(D_r^{d-1}).
\end{align*}
which means $D_r^{d-1}\setminus S(\Sigma^{\geqslant 2}(\phi)\cap C_{a,b}^d)\neq\varnothing$. It is equivalent to show that there exists some $z^*\in D_r^{d-1}$, such that $\text{int}\,(l_{z^*})\cap\Sigma^{\geqslant 2}(\phi)=\varnothing$. Thus, $l_{z^*}$ is a broken line from $a$ to $b$ we need.
\end{proof}

\begin{Cor}\label{cor:local_connected}
Assume that $\phi\in\emph{SCL}\,(M)$ and $d\geqslant 2$. Then for any $A\subseteq\Sigma^{\geqslant 2}(\phi)$, $M\setminus A$ is both path-connected and locally path-connected. In particular, for $k=1,2,\cdots d$, $\Sigma^{\leqslant k}(\phi)$ are both path-connected and locally path-connected subsets of $M$.
\end{Cor}

\begin{ex}
Actually, in some sense, $\Sigma^{\leqslant 1}(\phi)=\Sigma^0(\phi)\cup\Sigma^1(\phi)$ is the smallest subset with both path-connectedness and local path-connectedness.
\begin{enumerate}[$\bullet$]
	\item $\Sigma^0(\phi)$ does not have to be (locally) (path) connected.
	
	Consider the function
    \begin{align*}
	\phi_1(x_1,x_2,\cdots,x_d)=\begin{cases}
			(x_1+\frac{1}{2^{n-1}})(x_1+\frac{1}{2^{n}}),\,x_1\in[-\frac{1}{2^{n-1}},-\frac{1}{2^{n}}], n\in\N^*,\\
			0,\,\,\,\,\,\,\,\,\,\,\,\,\,\,\,\,\,\,\,\,\,\,\,\,\,\,\,\,\,\,\,\,\,\,\,\,\,\,\,\,\,\,\,\,\,\text{otherwise}.
		\end{cases}
    \end{align*}
    This is a semiconcave function in $\R^d$, and 
    \begin{align*}
	 \Sigma^0(\phi_1)=\left\{(x_1,x_2,\cdots,x_d)\in\R^d\big|x_1\neq-\frac{1}{2^{n-1}},n\in\N^*\right\},
    \end{align*}
    which is neither path-connected nor locally path-connected (even neither connected nor locally connected).
    \item $\Sigma^1(\phi)$ does not have to be (locally) (path) connected.
    
    Consider the function
    \begin{align*}
    	\phi_2(x_1,x_2,\cdots,x_d)=
    	\begin{cases}
    		(x_1+\frac{1}{2^{n-1}})(x_1+\frac{1}{2^{n}}),\,x_1\in[-\frac{1}{2^{n-1}},-\frac{1}{2^{n}}], n\in\N^*,\\
			-x_1,\,\,\,\,\,\,\,\,\,\,\,\,\,\,\,\,\,\,\,\,\,\,\,\,\,\,\,\,\,\,\,\,\,\,\,\,\,\,x_1\in[0,+\infty),\\
			0,\,\,\,\,\,\,\,\,\,\,\,\,\,\,\,\,\,\,\,\,\,\,\,\,\,\,\,\,\,\,\,\,\,\,\,\,\,\,\,\,\,\,\,\,\,\text{otherwise}.
    	\end{cases}
    \end{align*}
    This is a semiconcave function in $\R^d$, and
    \begin{align*}
    	\text{Sing}(\phi_2)=\Sigma^1(\phi_2)=\left\{(x_1,x_2,\cdots,x_d)\in\R^d\big| x_1=-\frac{1}{2^{n-1}},0,n\in\N^*\right\},
    \end{align*}
    which is neither path-connected nor locally path-connected (even neither connected nor locally connected).
\end{enumerate}
\end{ex}

It is still interesting to raise the following problems:
\begin{enumerate}[(a)]
	\item Are there some other consequences on the singular set of semiconcave functions by selecting a specific Hamiltonian in Arnaud's theorem?
	\item Is there an extension of Arnaud's theorem for semiconcave functions with general modulus?
	\item How about the singular set of semiconcave functions with general modulus? In fact, the conclusion in Corollary \ref{cor:local_connected} for the (local) path-connectedness is also true for the semiconcave function with general modulus by the discussion in \cite[Section 4.1]{Cannarsa_Sinestrari_book}.
\end{enumerate}

%\appendix

\subsection{Arnaud's theorem and Lasry-Lions regularization}

In this section we first discuss Proposition \ref{pro:Arnaud} and give a new and more geometric proof for time-dependent Hamiltonian. The following estimate is standard.% We give a proof in the appendix.

\begin{Lem}\label{lem:sup_max}
Let $\tau_2-\tau_1>0$ and let $\phi$ be a global Lipschitz function (with constant $\ell>0$). Then, there exists $\lambda_{\ell}>0$  satisfies the following property: if $x\in M$, $\tau_1\leqslant s<t\leqslant\tau_2$ and $y_{s,t,x}$ is a maximizer of the function $\phi(\cdot)-h(s,t,x,\cdot)$, then $y_{s,t,x}$ is contained in $B(x,\lambda_{\ell} (t-s))$. More precisely, 
\begin{align*}
	\lambda_{\ell}=c_1+\theta^*(\ell+1)+\max_{\substack{x\in M\\s\in[\tau_1,\tau_2]}}|L(s,x,0)|+c_0.
\end{align*}
\end{Lem}

\begin{Pro}\label{lem:uniqueness}
(See \cite[Appendix A]{Chen_Cheng_Zhang2018}) Under the same assumption as Lemma \ref{lem:sup_max} and $\lambda_{\ell}$ is given by Lemma \ref{lem:sup_max}. Then, there exists $t_{\phi}>0$ such that for any $s\in[\tau_1,\tau_2)$ and $x\in M$, %$\phi\in\emph{SCL}\,(M)$ and $x\in M$,
\begin{enumerate}[\rm (1)]   
    \item the functions $(t,y)\mapsto h(s,t,x,y)$ and  $(t,y)\mapsto h(s,t,y,x)$ are both locally semiconcave and locally semiconvex (thus $C^{1,1}_{\emph{loc}}$) on
        \begin{align*}
        	S(s,x,\lambda_\ell,t_\phi)=\{(t,y)\in\R\times M: 0<t-s\leqslant\max\{s+t_\phi,\tau_2\}, d(y,x)\leqslant \lambda_\ell(t-s)\}.
        \end{align*}
        For all $(t,y)\in S(s,x,\lambda_\ell,t_\phi)$,
        \begin{align*}
        	 D_yh(s,t,x,y)&=L_v(t,\xi(t),\dot\xi(t)),\\
        	 D_xh(s,t,x,y)&=-L_v(s,\xi(s),\dot\xi(s)),
        \end{align*}
        where $\xi\in\Gamma_{x,y}^{s,t}$ is the unique minimizer of $h(s,t,x,y)$;
	\item the function $\phi(\cdot)-h(s,t,x,\cdot)$ is strictly concave in $B(x,\lambda_{\ell} (t-s))$ for all $s<t\leqslant\max\{s+t_\phi,\tau_2\}$;
	\item the function $\phi(\cdot)-h(s,t,x,\cdot)$ admits a unique maximizer on $B(x,\lambda_{\ell} (t-s))$ (so on $M$) for all $s<t\leqslant\max\{s+t_\phi,\tau_2\}$.
\end{enumerate}
\end{Pro}

Fix $t_1<t_2$, since the local nature of the following discussion, we can consider the characteristic system on $[t_1,t_2]$ as
\begin{equation}\label{eq:Lie}\tag{C}
	\begin{cases}
			\dot X(t)=H_p(t,X(t),P(t)),\\
			\dot P(t)=-H_x(t,X(t),P(t)),
	\end{cases}
	\qquad\text{and}\ \qquad
	\begin{cases}
			X(t_2)=x,\\
			P(t_2)=p.
		\end{cases}
\end{equation}
We denote the solution of \eqref{eq:Lie} by
\begin{align*}
(X(\cdot),P(\cdot)):=(X(\cdot;t_1,t_2,x,p),P(\cdot;t_1,t_2,x,p)).
\end{align*}	

\begin{Lem}\label{lem:diff}
Fix $\tau_2>\tau_1$. For any $R>0$ there exists $m_R, t_R>0$ such that, if $x\in M$, $\tau_1\leqslant t_1<t_2\leqslant\tau_2$ and $0<t_2-t_1\leqslant t_R$, then the map
\begin{align*}
	\Phi_{x,t_1,t_2}:	 B(0,R)&\to B(x,m_Rt_R)\\
	p&\mapsto X(t_1;t_1,t_2,x,p)
\end{align*}
is a $C^2$-diffeomorphism.	
\end{Lem}

\begin{proof}
Consider the variational equation of \eqref{eq:Lie} with respect to $p$
\begin{equation}\label{eq:var_p}
	\begin{cases}
			\dot X_p(t)=H_{px}(t,X(t),P(t))X_p(t)+H_{pp}(t,X(t),P(t))P_p(t),\\
			\dot P_p(t)=-H_{xx}(t,X(t),P(t))X_p(t)-H_{xp}(t,X(t),P(t))P_p(t),\\
		\end{cases}
\end{equation}
with terminal condition $X_p(t_2)=0$, $P_p(t_2)=I$.  It is easy to check that $D\Phi_{x,t_1,t_2}(p)= X_p(t_1).$

From \eqref{eq:var_p}, we can get $\dot X_p(t_2)=H_{pp}(t_2,x,p)$. Since $H$ is strictly convex with respect to component $p$, we conclude that for any $R>0$, there exists $c_R>0$, such that $H_{pp}(t_2,x,p)>c_RI$ when $|p|_x\leqslant R$. Moreover, by the continuity of $\dot X_p(\cdot)$, there exists $t_R>0$ such that, if $0<t_2-t_1<t_R$,
\begin{align*}
	\dot X_p(r)>\frac{c_R}{2}I,\qquad\forall r\in[t_1,t_2].
\end{align*}
Thus, for $|p|_x\leqslant R$ and $0<t_2-t_1\leqslant t_R$ we obtain
\begin{align*}
	X_p(t_1;t_1,t_2,x,p)=X_p(t_2;t_1,t_2,x,p)-\int^{t_2}_{t_1}\dot X_p(r)\ dr<-\frac{c_R(t_2-t_1)}{2}I,
\end{align*}
and our conclusion follows from the inverse mapping theorem.
\end{proof}

\begin{The}\label{thm:C11}
Suppose $\phi\in\emph{SCL}\,(M)$ and $\tau_2-\tau_1>0$. Then there exist $t_\phi'$ and $t_\phi$ with $t_\phi^\prime\geqslant t_\phi>0$ independent of $x$ such that the following properties hold.
\begin{enumerate}[\rm (1)]
	\item If $\tau_1\leqslant t_1<t_2\leqslant\tau_2$ and $t_2-t_1\leqslant t_\phi^\prime$, then $\Phi_{H}^{t_1,t_2}(\text{\rm graph}\,(D^*\breve{T}_{t_1}^{t_2}\phi))\subseteq\text{\rm graph}\,(D^+\phi)$.
	\item If $\tau_1\leqslant t_1<t_2\leqslant\tau_2$ and $t_2-t_1\leqslant t_\phi$, then there exists a family $\mathscr F_\phi\subseteq C^2(M)$ with $\|f\|_{C^2}$ uniformly bounded such that
	\begin{equation}\label{eq:min1}
	    \breve T_{t_1}^{t_2}\phi=\inf_{f\in\mathscr F_\phi}\breve T_{t_1}^{t_2}f.
	\end{equation}
	\item If $\tau_1\leqslant t_1<t_2\leqslant\tau_2$ and $t_2-t_1\leqslant t_\phi$, then $\breve T_{t_1}^{t_2}\phi\in C^{1,1}(M)$.
	\item The inclusion in \mbox{\rm (1)} is indeed an equality. That is, If $\tau_1\leqslant t_1<t_2\leqslant\tau_2$ and $t_2-t_1\leqslant t_\phi$, then
	\begin{equation}\label{eq:graph}
	    \Phi_{H}^{t_1,t_2}(\text{\rm graph}\,(D\breve{T}_{t_1}^{t_2}\phi))=\text{\rm graph}\,(D^+\phi).
	\end{equation}
\end{enumerate}
\end{The}

\begin{proof}

Fix $x\in M$ and suppose $0<\tau_2-\tau_1\leqslant1$ and $\tau_1\leqslant t_1<t_2\leqslant\tau_2$, then Lemma \ref{lem:sup_max}, Proposition \ref{lem:uniqueness} and Lemma \ref{lem:diff} imply there exist $\lambda>0$ and $t'_\phi\geqslant t_\phi>0$ such that
\begin{enumerate}[--]
	\item the function $h(t_1,t_2,x,\cdot)$ is of class $C^{1,1}_{\rm loc}(B(x,\lambda (t_2-t_1)))$ for all $0<t_2-t_1\leqslant t'_\phi$;
	\item the function $\phi(\cdot)-h(t_1,t_2,x,\cdot)$ is strictly concave in $B(x,\lambda (t_2-t_1))$ for all $0<t_2-t_1\leqslant t_\phi$;
	\item $\Phi_{t_1,t_2}$ is a $C^2$-diffeomorphism when $0<t_2-t_1\leqslant t_\phi$.
\end{enumerate}

Now we suppose $t_2-t_1\in(0,t'_\phi]$. Using the notation in Proposition \ref{pro:marginal}, then
\begin{align*}
	M(x)=\arg\max\{\phi(y)-h(t_1,t_2,x,y): y\in M\},\quad Y(x)=\{-D_xh(t_1,t_2,x,y): y\in M(x)\},
\end{align*}
and $D^*\breve{T}_{t_1}^{t_2}\phi(x)\subseteq Y(x)$. To prove (1), note that for any $p\in D^*\breve{T}_{t_1}^{t_2}\phi(x)$ there exists $y\in M(x)$ such that
\begin{align*}
	\breve{T}_{t_1}^{t_2}\phi(x)=\phi(y)-h(t_1,t_2,x,y)=\phi(\gamma(t_2))-\int_{t_1}^{t_2}L(s,\gamma,\dot\gamma)\,ds,
\end{align*}
where $\gamma\in\Gamma^{t_1,t_2}_{x,y}$ is the unique minimal curve for $h(t_1,t_2,x,y)$, and
\begin{align*}
	p=-D_xh(t_1,t_2.x,y)=L_v(t_1,\gamma(t_1),\dot{\gamma}(t_1))=p(t_1)
\end{align*}
by Proposition \ref{lem:uniqueness}, where $p(s):=L_v(s,\gamma(s),\dot{\gamma}(s))$ is the dual arc of $\gamma$. On the other hand,  by Fermat's rule, 
\begin{align*}
	0\in D_y^+(\phi(y)-h(t_1,t_2,x,y))=D_y^+\phi(y)-D_y h(t_1,t_2,x,y).
\end{align*}
By Proposition \ref{lem:uniqueness} again, the relation above reads
\begin{equation}\label{eq:Fermat_rule}
	p(t_2)=L_v(t_2,\gamma(t_2),\dot\gamma(t_2))\in D^+\phi(y).
\end{equation}
This leads to the inclusion in (1).

Now, suppose $t_2-t_1\in(0,t_\phi]$. We turn to the proof of (2). Since $\phi\in \text{SCL}\,(M)$, for any $x\in M$ and $p\in D^+\phi(x)$ one can find a $C^2$ function $f$ touching $\phi$ from above at $x$. That is $f(x)=\phi(x)$ and $p=Df(x)$. We denote by $\mathscr F_\phi$ the set of such $C^2$ functions. Without loss of generality we assume that $\mathscr F_\phi$ is bounded with the $C^2$-norm since $\phi$ is (globally) semiconcave with linear modulus. Then it follows
\begin{align*}
	\phi(x)=\min_{f\in\mathscr F_\phi}f(x).
\end{align*}
It is not difficult to see that $\breve T_{t_1}^{t_2}\phi\leqslant \breve T_{t_1}^{t_2}f$  for any $f\in\mathscr F_\phi$. Thus $\breve T_{t_1}^{t_2}\phi\leqslant\min_{f\in\mathscr F_\phi}\breve T_{t_1}^{t_2}f$. On the other hand, one can find $f\in\mathscr F_\phi$ touching $\phi$ from above at $y$ in the proof of (1) and we observe
\begin{align*}
	\gamma(s)=\pi_x\Phi_{H}^{t_1,s}(y,Df(y)),\qquad s\in[t_1,t_2].
\end{align*}
It follows
\begin{align*}
	\breve T_{t_1}^{t_2}f(x)&=f(\gamma(t_2))-\int_{t_1}^{t_2}L(s,\gamma(s),\dot \gamma(s))\ ds\\
&=\phi(y)-h(t_1,t_2,x,y)\\
&=\breve T_{t_1}^{t_2}\phi(x),
\end{align*}
It follows $\breve{T}_{t_1}^{t_2}\phi(x)=\breve T_{t_1}^{t_2}f(x)\geqslant\inf_{f\in\mathscr F_\phi}\breve T_{t_1}^{t_2}f
(x)$. This completes the proof of (2). 

The conclusion in (3) is a direct consequence of (2). Indeed $\breve{T}_{t_1}^{t_2}\phi$ is semiconcave by \eqref{eq:min1}.  As the supreme of the family of semiconvex functions $\phi(y)-h(t_1,t_2,\cdot,y)$ it is also semiconvex.

To obtain \eqref{eq:graph}, it is enough to prove the converse direction, i.e., $\Phi_{H}^{t_1,t_2}(\text{\rm graph}\,(D\breve{T}_{t_1}^{t_2}\phi))\supseteq\text{\rm graph}\,(D^+\phi)$. Choose $R>0$ such that $D^+\phi(y)\subseteq B(0,R)$ and let $t_R$ be given in Lemma \ref{lem:diff}. For any $t_2-t_1\in(0,t_R]$, $y\in M$ and $p\in D^+\phi(y)$, we denote $\xi(s)=X(s;t_1,t_2,y,p)$, $s\in[t_1,t_2]$, $x=\xi(t_1)$. By Lemma \ref{lem:diff}, $\xi$ is the unique minimizer for $h(t_1,t_2,x,y)$. Combining Proposition \ref{lem:uniqueness} (1), we conclude
\begin{align*}
	p=L_v(t_2,\xi(t_2),\dot{\xi}(t_2))=D_yh(t_1,t_2,x,y),
\end{align*}
when $t_R\ll1$. This implies
\begin{align*}
	D_yh(t_1,t_2,x,y)\in D^+\phi(y).
\end{align*}
In other words, $y$ is the unique maximizer of the function $\phi(\cdot)-h(t_1,t_2,x,\cdot)$, by Proposition \ref{lem:uniqueness} (2) and (3). Therefore
\begin{align*}
	\Phi_{H}^{t_1,t_2}((x,D\breve{T}_{t_1}^{t_2}\phi(x))=(y,p).
\end{align*}
This completes the proof of (4).
\end{proof}

\bibliographystyle{plain}
\bibliography{mybib}

\begin{thebibliography}{10}

\bibitem{Albano_Cannarsa1999}
Paolo Albano and Piermarco Cannarsa.
\newblock Structural properties of singularities of semiconcave functions.
\newblock {\em Ann. Scuola Norm. Sup. Pisa Cl. Sci. (4)}, 28(4):719--740, 1999.

\bibitem{Albano_Cannarsa2000}
Paolo Albano and Piermarco Cannarsa.
\newblock Propagation of singularities for concave solutions of
  {H}amilton-{J}acobi equations.
\newblock In {\em International {C}onference on {D}ifferential {E}quations,
  {V}ol. 1, 2 ({B}erlin, 1999)}, pages 583--588. World Sci. Publ., River Edge,
  NJ, 2000.

\bibitem{Alberti_Ambrosio_Cannarsa1992}
Giovanni Alberti, Luigi Ambrosio, and Piermarco Cannarsa.
\newblock On the singularities of convex functions.
\newblock {\em Manuscripta Math.}, 76(3-4):421--435, 1992.

\bibitem{Arnaud2011}
M.-C. Arnaud.
\newblock Pseudographs and the {L}ax-{O}leinik semi-group: a geometric and
  dynamical interpretation.
\newblock {\em Nonlinearity}, 24(1):71--78, 2011.

\bibitem{Cannarsa_Cheng3}
Piermarco Cannarsa and Wei Cheng.
\newblock Generalized characteristics and {L}ax-{O}leinik operators: global
  theory.
\newblock {\em Calc. Var. Partial Differential Equations}, 56(5):Art. 125, 31,
  2017.

\bibitem{Cannarsa_Cheng_Fathi2017}
Piermarco Cannarsa, Wei Cheng, and Albert Fathi.
\newblock On the topology of the set of singularities of a solution to the
  {H}amilton-{J}acobi equation.
\newblock {\em C. R. Math. Acad. Sci. Paris}, 355(2):176--180, 2017.

\bibitem{Cannarsa_Cheng_Fathi2021}
Piermarco Cannarsa, Wei Cheng, and Albert Fathi.
\newblock Singularities of solutions of time dependent {H}amilton-{J}acobi
  equations. {A}pplications to {R}iemannian geometry.
\newblock {\em Publ. Math. Inst. Hautes \'{E}tudes Sci.}, 133(1):327--366,
  2021.

\bibitem{Cannarsa_Sinestrari_book}
Piermarco Cannarsa and Carlo Sinestrari.
\newblock {\em Semiconcave functions, {H}amilton-{J}acobi equations, and
  optimal control}, volume~58 of {\em Progress in Nonlinear Differential
  Equations and their Applications}.
\newblock Birkh{\"a}user Boston, Inc., Boston, MA, 2004.

\bibitem{Cannarsa_Yu2009}
Piermarco Cannarsa and Yifeng Yu.
\newblock Singular dynamics for semiconcave functions.
\newblock {\em J. Eur. Math. Soc. (JEMS)}, 11(5):999--1024, 2009.

\bibitem{Chen_Cheng_Zhang2018}
Cui Chen, Wei Cheng, and Qi~Zhang.
\newblock Lasry-{L}ions approximations for discounted {H}amilton-{J}acobi
  equations.
\newblock {\em J. Differential Equations}, 265(2):719--732, 2018.

\bibitem{Edgar_book2008}
Gerald Edgar.
\newblock {\em Measure, topology, and fractal geometry}.
\newblock Undergraduate Texts in Mathematics. Springer, New York, second
  edition, 2008.

\bibitem{Engelking_book1978}
Ryszard Engelking.
\newblock {\em Dimension theory}, volume~19 of {\em North-Holland Mathematical
  Library}.
\newblock North-Holland Publishing Co., Amsterdam-Oxford-New York; PWN---Polish
  Scientific Publishers, Warsaw, 1978.
\newblock Translated from the Polish and revised by the author.

\bibitem{Fathi_Figalli2010}
Albert Fathi and Alessio Figalli.
\newblock Optimal transportation on non-compact manifolds.
\newblock {\em Israel J. Math.}, 175:1--59, 2010.

\bibitem{Federer_book1969}
Herbert Federer.
\newblock {\em Geometric measure theory}.
\newblock Die Grundlehren der mathematischen Wissenschaften, Band 153.
  Springer-Verlag New York, Inc., New York, 1969.

\bibitem{Hiriart-Urruty_Lemarechal_book2001}
Jean-Baptiste Hiriart-Urruty and Claude Lemar{\'e}chal.
\newblock {\em Fundamentals of convex analysis}.
\newblock Grundlehren Text Editions. Springer-Verlag, Berlin, 2001.
\newblock Abridged version of {{\i}t Convex analysis and minimization
  algorithms. I} [Springer, Berlin, 1993; MR1261420 (95m:90001)] and {{\i}t II}
  [ibid.; MR1295240 (95m:90002)].

\bibitem{Hurewicz_Wallman1941}
Witold Hurewicz and Henry Wallman.
\newblock {\em Dimension {T}heory}.
\newblock Princeton Mathematical Series, vol. 4. Princeton University Press,
  Princeton, N. J., 1941.

\bibitem{Zajicek1979}
Lud\v{e}k Zaj\'{\i}\v{c}ek.
\newblock On the differentiation of convex functions in finite and infinite
  dimensional spaces.
\newblock {\em Czechoslovak Math. J.}, 29(104)(3):340--348, 1979.

\end{thebibliography}

\end{document}